\documentclass{article}[12pt]
\usepackage{amsfonts,amssymb}
\newtheorem{theorem}{Theorem}

\newtheorem{prop}{Proposition}

\def\PSL{\mbox{\rm{PSL}}}
\def\Ker{\mbox{\rm{Ker}}}

\begin{document}

\title{Kleinian groups with discrete length spectrum}
\author{Richard D. Canary\thanks{Research supported in part by
NSF grants 0203698 and 0504791}\\ Christopher J. Leininger\thanks{Research
supported in part by an N.S.F. postdoctoral fellowship}}

\maketitle

Given a Kleinian group
$\Gamma \subset \PSL_{2}({\mathbb C})\cong {\rm Isom}_+({\mathbb H}^3)$,
we define the {\em real length spectrum of $\Gamma$ without multiplicities}
to be the set of all real translation lengths of elements of
$\Gamma$ acting on ${\mathbb H}^{3}$
$${\frak L}(\Gamma) = \{ \, l(x) \, | \, x \in \Gamma \, \} \subset {\mathbb R}$$
Equivalently, these are the lengths of closed geodesics in
the quotient orbifold.

In this paper we characterize finitely generated, torsion-free
Kleinian groups with discrete real length spectrum.  In the process,
we obtain a sharper version of the Covering Theorem which may
be of independent interest.

\begin{theorem} \label{main1}
Suppose that $\Gamma$ is a finitely generated, torsion-free
Kleinian group.
Then the real length spectrum ${\frak L}(\Gamma)$ is discrete
if and only if either
\begin{enumerate}
\item $\Gamma$ is geometrically finite,
\item there exists a hyperbolic 3-manifold $M=\mathbb{H}^3/G$
which fibers over the circle and $\Gamma$ is a fiber subgroup of $G$, or
\item there exists a hyperbolic 3-manifold $M=\mathbb{H}^3/G$ which fibers over
\newline
$S^{1}/\langle z \mapsto \overline{z} \rangle$ and $\Gamma$ is a singular
fiber subgroup of $G$.
\end{enumerate}
\end{theorem}

We recall that if a hyperbolic 3-manifold $M=\mathbb{H}^3/G$
fibers over a 1-orbifold $Q$, then a fiber subgroup of $G$
is the fundamental group of the pre-image of a regular
point of $Q$, while a singular fiber subgroup is the fundamental group
of the pre-image of a singular point of $Q$. (Notice that 
one only really obtains a conjugacy class of subgroups.)

\medskip
\noindent
{\bf Proof:}
We refer the reader to \cite{Ca1} and \cite{Ca2} for details and terminology.
Let $N = {\mathbb H}^{3}/\Gamma$, let $\epsilon > 0$ be less than the
$3$-dimensional Margulis constant, and let $N_{\epsilon}^0$ denote
the complement of the $\epsilon$-cuspidal thin part of $N$.
By the Tameness Theorem (see Agol \cite{Ag} or  Calegari-Gabai \cite{CG}),
$N$ is topologically tame.

It is well-known that a geometrically finite Kleinian group has
discrete real length spectrum.
We sketch a brief proof for the reader's convenience. If $\Gamma$
is geometrically finite, then there exists a compact
subset $K$ of $N=\mathbb{H}^3/\Gamma$ such that any closed geodesic 
intersects $K$. (One may, for example,  take $K$ to be the intersection
of the convex core of $N$ with $N^0_\epsilon$.)
If ${\frak L}(\Gamma)$ is indiscrete, then there exists an infinite
collection of distinct closed geodesics of length at most $L$ (for
some $L>0$) all of which intersect $K$. If this were the case,
some sequence of closed geodesics in $N$ would have to accumulate
at a closed geodesic, which is impossible.
Therefore, ${\frak L}(\Gamma)$ is discrete.

Since a fiber subgroup, as described in (2), or singular fiber subgroup,
as described in (3), is a subgroup of a geometrically finite Kleinian group,
they must also have discrete real length spectrum.

We now assume that $\Gamma$ is a geometrically infinite Kleinian
group with discrete real length spectrum and prove that $\Gamma$ is
either a fiber subgroup or a singular fiber subgroup
of the fundamental group of
a finite volume hyperbolic 3-manifold which fibers over a 1-orbifold.
Let $U \cong \overline{S} \times [0,\infty) \subset N_{\epsilon}^0$
be a neighborhood of a geometrically infinite end of
$N^0_\epsilon$ and let $U_p \cong S \times [0,\infty) \subset N$
be its parabolic extension.
By the main result of \cite{Ca1}, there exists a sequence of useful simplicial
hyperbolic surfaces $f_{n}:X_{n} \rightarrow U_p$, each properly homotopic
in $U_p$ to $S \times \{ 0 \}$, which exits the end.

As the length spectrum of $\Gamma$ is discrete, there is an uniform positive
lower bound to the length of any closed geodesic in $N$. Therefore, there
is an uniform positive lower bound on the
length of any closed geodesic on $X_{n}$ whose image
under $f_n$ is essential and not an accidental parabolic (i.e. not
homotopic into a cusp of $N$.).
On the other hand, Lemma 7.1 in \cite{Ca1} guarantees that
given any $A > 0$, there exists $n_A > 0$ so that for all
$n \geq n_A$, any closed geodesics in $X_n$ which is inessential
or accidentally parabolic under $f_n$ must have length at least
$A$ on $X_n$.  Therefore, there is an uniform lower bound to the
length of a closed geodesic on any $X_n$.

Lemma 7.1 in \cite{Ca2} then gives that there exists $C>0$ and,
for all $n$, a point $x_n$ in the $\epsilon_0$-thick part of $X_n$
and a minimal generating set $\{ a_1^n,\ldots,a_m^n \}$ for $\pi_1(X_n,x_n)$
all of which are represented by curves (based at $x_n$) of length
at most $C$. By passing to a subsequence, we may assume that all
the minimal generating sets are topologically equivalent.
Therefore, there exist markings $g_{n}:S \rightarrow X_{n}$
so that any chosen finite set of curves on $S$ have representatives with
uniformly bounded length on $X_{n}$.
Composing, we obtain a sequence of representations
$$\rho_n=(f_{n} \circ g_{n})_{*}:\pi_{1}(S) \rightarrow \Gamma \subset
\PSL_{2}({\mathbb C})$$

A result of Kim \cite{Ki} gives that there exists a
finite set $\{s_{1},...,s_{k}\}$ of elements of $\pi_{1}(S)$ such that
an irreducible representation of $\pi_{1}(S)$
is determined up to conjugacy by their lengths.
Since the lengths of $s_{1},...,s_{k}$ with respect to
$g_{n}:S \rightarrow X_{n}$ are bounded, so are the set of lengths
$\rho_n(s_{i})$. Since ${\frak L}(\Gamma)$ is discrete,  we may pass to
a subsequence, again called $\{\rho_n\}$, so that $l(\rho_n(s_i))$
is constant for all $i$.  Therefore, all representations in this sequence
are conjugate to $\rho_1$.
Let $h_{n} \in \PSL_{2}({\mathbb C})$ be the conjugating element, so that
$\rho_n= h_{n} \rho_1h_{n}^{-1}$.

We now observe that $U_p$ is incompressible (i.e. that $\rho_1$ is injective).
First note that, since all the $\rho_n$ are conjugate, if
$s \in \Ker(\rho_1)$ is non-trivial,
then $s \in \Ker(\rho_n)$ for all $n$, and that
the length of $g_{n}(s)$ in $X_{n}$ is uniformly bounded.
However, since $f_n(g_n(s))$ is non-essential, $\{ l_{X_n}(g_n(s))\}$ converges
to $\infty$, which is a contradiction.

Finally, consider the cover $\widetilde N_U$ of $N$
corresponding to $\pi_{1}(U)$.
The conjugating elements $h_{n}$ descend to isometries
$\hat h_n$ of $\widetilde N_U$.
Let $\widetilde f_n:X_n\to \widetilde N_U$ denote the lift of $f_n$
to $\widetilde N_U$ and let $\widetilde x_n =\widetilde f_n (x_n)$.
Note that there exists $D>0$ so that if
$x\in\widetilde N_U$ and there exist representatives of
$\{\widetilde f_1(a_1^1), \ldots, \widetilde f_1(a_m^1)\}$ based at $x$ of length at most $C$, then
$d(x,\widetilde x_1)\le D$. It then follows that
$d(\widetilde x_n,h_n(\widetilde x_1))\le D$ for all $n$.
Since $\{f_n(X_n)\}$ exits $U$, one sees that $\hat h_n(U)\cap U$ is non-empty
for all large enough $n$, and that there are infinitely many distinct
$\hat h_n$.

Let $N_1$ be the quotient of $\widetilde N_U$ by its full group
of orientation-preserving isometries. We have established that
the covering map $\pi_U:\widetilde N_U\to N_1$ is
infinite-to-one on $U$.

We now recall the  Covering Theorem (see \cite[Theorem 9.2.2]{T}, \cite{Ca2}
and Agol \cite{Ag}):

\medskip\noindent
{\bf Covering Theorem:} {\em Suppose that $\widehat M$ is a hyperbolic
3-manifold with finitely generated fundamental group and $M$
is a hyperbolic 3-orbifold. If $p:\widehat M\to M$ is
an orbifold cover  which is infinite-to-one on a neighborhood $U$ of
a geometrically infinite end of $\widehat M^0_\epsilon$, then 
$M$  has finite volume and has a finite manifold cover $M'$
which fibers over the circle such that either

a) $\widehat M$  is the cover of $M'$ associated to the fiber, or

b) $\widehat M_0^\epsilon$ is a twisted ${\bf R}$-bundle and
$\widehat M$ is  double covered by the manifold $\widehat M'$
which is the cover of $M'$ associated to the fiber.}

\medskip

The Covering Theorem immediately implies that $\Gamma_U$ is the fiber
group associated to a 3-manifold which fibers over the circle.
If $\Gamma=\Gamma_U$, we're done.
If not, we apply the Covering Theorem again to the map
$p_U:\widetilde N_U\to N$ to conclude that $p_U$ is finite-to-one.
Hempel's Finite Index Theorem \cite[Theorem 10.5]{Hempel} implies
that $N_0^\epsilon$ is a twisted $\bf R$-bundle and that the cover $p_U$ is
two-to-one.

The following addendum to the Covering Theorem allows
us to conclude that $\Gamma$ is a singular fiber subgroup,
in the case that $N^0_\epsilon$ is a twisted $\bf R$-bundle.

\begin{prop}
\label{addendum}
Suppose that $\widehat M$ is a geometrically infinite
hyperbolic 3-manifold such that $\widehat M_0^\epsilon$
is a twisted $\bf R$-bundle over a compact
surface, and let $p':\widehat M' \to \widehat M$ be a two-fold cover
such that $\widehat M_0^\epsilon$ an untwisted $\bf R$-bundle. 

\medskip\noindent
1) If there is a finite volume hyperbolic 3-manifold $M''$ that fibers over the circle with covering $q':\widehat M' \to M''$ corresponding to the fiber, then $\widehat M$ covers a finite volume hyperbolic orbifold $M$.

\medskip\noindent
2) If $p:\widehat M \to M$ is a cover of a finite volume hyperbolic
orbifold $M$, then $M$ has a finite manifold cover $M'=\mathbb{H}^3/G'$
which fibers over the orbifold $S^{1}/\langle z \mapsto \overline{z} \rangle$
and $\widehat M$ is the cover associated to a singular fiber subgroup of $G'$.
\end{prop}

\newpage

We now turn to the proof of Proposition \ref{addendum}:
For the first statement, we note that the covers
$p':\widehat M' \to \widehat M$ and $q':\widehat M' \to M''$
are regular covers with covering groups generated by an involution
$a \in {\rm Isom }_+(\widehat M')$ and infinite order isometry
$b \in {\rm Isom}_+(\widehat M')$, respectively.
Since ${\rm Isom }_+(\widehat M')$ acts properly discontinuously
on $\widehat M'$ and $\langle b \rangle$ has quotient $M''$ of
finite volume, so the covering
$\hat p: \widehat M' \to M = \widehat M'/{\rm Isom }_+(\widehat M')$
is onto the finite volume orbifold $M$.
Since $\widehat M$ is the quotient of $\widehat M'$ by the subgroup
$\langle a \rangle < {\rm Isom}_+(\widehat M')$, it follows that there
is a covering $p:\widehat M \to M$ so that $\hat p = p \circ p'$,
verifying the first statement.

For the second assertion, define $\hat p = p \circ p'$. 
The Covering Theorem, applied to $\hat p:\widehat M'\to M$,
implies that $M$ has a finite manifold cover $M''$ which fibers over the circle
so that $\widehat M'$ is the cover associated to the fiber.
Moreover, examining the proof of the Covering Theorem,
we see that if $q:M''\to M$ and $q':\widehat M'\to M''$ are the
associated covering maps, we may assume that $q\circ q'=\hat p=p\circ p'$. 

As above, the covers $p':\widehat M' \to \widehat M$ and
$q':\widehat M' \to M''$ are regular with covering group of order two,
generated by $a \in {\rm Isom}_+(\widehat M')$, and infinite cyclic
covering group, generated by $b \in {\rm Isom}_+(\widehat M')$,
respectively.  Since $\langle b \rangle$ has finite volume quotient
$M''$ and ${\rm Isom}_+(\widehat M')$ acts properly discontinuously,
it follows $\langle b\rangle$ has finite index in ${\rm Isom}_+(\widehat M')$
and hence in $\langle a,b\rangle$. Therefore, there exists $k>0$ such
that $\langle b^k\rangle$ is normal in $\langle a, b\rangle$.

Since $\langle b^k\rangle$ is normal in $\langle a,b\rangle $, we see
that either $ab^ka=b^{k}$ or $ab^ka=b^{-k}$. We now see that it
must be the case that $ab^ka=b^{-k}$. Let $T$ be a connected
separating surface in $(\widehat M')^0_\epsilon$ which is preserved
by $a$, i.e. $a(T)=T$. If $n>0$
is chosen large enough then $b^{nk}(T)$ does not intersect $T$. Since
$a$ interchanges the two components of $(N_U)^0_\epsilon - T$,
we see immediately that $ab^{nk}a(T)$ and $b^{nk}(T)$ lie on opposite
sides of $T$. Therefore, $ab^{nk}a\ne b^{nk}$, which implies
that $ab^ka\ne b^k$. It follows that $ab^ka=b^{-k}$ and
that $\langle a,b^k \rangle \cong \mathbb{Z}_2*\mathbb{Z}_2$ is generated
by $a$ and $ab^k$.

If $n$ is even, then $ab^{nk}$ is conjugate in $\langle a,b\rangle$ to $a$.
Since $a$ acts freely, so does $ab^{nk}$. Therefore,
$\langle a,b^{2k}\rangle\cong \mathbb{Z}_2*\mathbb{Z}_2$ acts freely on $M'$.
Let $M'$ be the
quotient of $\widehat M'$ by the group $\langle a,b^{2k}\rangle$.
The main result of Hempel-Jaco \cite{HJ}
then implies that $M'$ fibers over the orbifold
$S^{1}/\langle z \mapsto \overline{z} \rangle$ and that $\widehat M'$
is the cover associated to a singular fiber subgroup.
(Alternatively, one may assume that  $k$ is  large enough that $b^{k}(T)$
does not intersect $T$ and prove directly that the region between $T$
and $b^{k}(T)$ is a fundamental region for the group $\langle a,b^{2k}\rangle.$
One may then fairly explicitly check that the quotient
$\widehat M'/ \langle a,b^{2k}\rangle$ fibers over the 
orbifold $S^{1}/\langle z \mapsto \overline{z} \rangle$ and that
$\widehat M'$ is a cover associated to a singular fiber subgroup.)
One may then check that $\hat p:\widehat M'\to M$ descends to
a covering map $\hat p':M'\to M$. This completes the proof of
Proposition \ref{addendum}, which in turn completes the proof
of our Main Theorem.

\bigskip

If we combine Proposition \ref{addendum} with our earlier
statement of the Covering Theorem, we get the following
slightly sharper version.

\newpage

\noindent
{\bf Covering Theorem(sharper version):} {\em Suppose that $\widehat M$ is a hyperbolic
3-manifold with finitely generated fundamental group and $M$
is a hyperbolic 3-orbifold. If $p:\widehat M \to M$ is
an orbifold cover  which is infinite-to-one on a neighborhood $U$ of
a geometrically infinite end of $\widehat M^0_\epsilon$, then 
$M$ has finite volume and has a finite manifold cover $M'=\mathbb{H}^3/G'$
such that either
\begin{enumerate}
\item
$M'$ fibers over the circle and $\widehat M$ is the cover associated
to a fiber subgroup of $G'$, or
\item
$M'$ fibers over the orbifold
$S^{1}/\langle z \mapsto \overline{z} \rangle$ and
$\widehat M$  is the cover of $M'$ associated to a singular fiber subgroup
of $G'$.
\end{enumerate}
}


\end{document}